\documentclass[reqno,10pt]{amsart}

\usepackage[left=3cm,right=3cm,top=3cm,bottom=3cm]{geometry}
\usepackage{amssymb,bm}
\usepackage{mathrsfs}
\usepackage{xcolor}
\usepackage{enumitem}
\usepackage{appendix}
\usepackage{diffcoeff}
\usepackage{amsmath,cases}

\usepackage{graphicx}
\usepackage{subfig}

\newtheorem{theorem}{Theorem}[section]

\theoremstyle{definition}

\theoremstyle{remark}
\newtheorem{remark}[theorem]{Remark}

\numberwithin{equation}{section}

\begin{document}

\title[Statistical conservation laws]{Statistical conservation laws for scalar model problems: Hierarchical evolution equations}



\author[Qian Huang]{Qian Huang}
\address{Institute of Applied Analysis and Numerical Simulation, University of Stuttgart, 70569 Stuttgart, Germany}
\email{qian.huang@mathematik.uni-stuttgart.de; hqqh91@qq.com}

\author[Christian Rohde]{Christian Rohde*}
\address{Institute of Applied Analysis and Numerical Simulation, University of Stuttgart, 70569 Stuttgart, Germany}
\email{christian.rohde@mathematik.uni-stuttgart.de}
\thanks{* Corresponding author}



\date{}

\dedicatory{}

\begin{abstract}
  The probability density functions (PDFs) for the solution of the incompressible Navier-Stokes equation can be represented by a hierarchy of linear equations. This article develops new hierarchical evolution equations for PDFs of a scalar conservation law with random initial data as a model problem. Two frameworks are developed, including multi-point PDFs and single-point higher-order derivative PDFs. These hierarchies capture statistical correlations and guide closure strategies.
\end{abstract}

\maketitle

\section{Introduction}

Consider the initial-value problem for $u = u(t,x, \xi(\omega)) \in \mathbb R$ with random initial data:
\begin{equation} \label{eq:scalar}
\begin{split}
  & u_t + \nabla_x\cdot g(u) = \epsilon \Delta_x u, \\
  & u(0,x,\xi(\omega)) = u_0(x,\xi(\omega)),
\end{split}
\quad t>0, x\in D \subset \mathbb R^d, \omega\in \Omega.
\end{equation}
Here $\xi: \Omega\to \mathbb R$ is a scalar absolutely-continuous random variable with a density function $p_\xi$, $g=g(u)\in\mathbb R^d$ is smooth, and $\epsilon>0$ is constant. We are interested in the statistical conservation laws for this scalar model problem.

An important reason to study (\ref{eq:scalar}) is that it mimics, at the scalar level, the interplay of nonlinear hyperbolic transport and viscosity operators that plays a key role for understanding turbulent flows. To tackle the statistical nature of turbulence, the efforts based on probabilistic approaches to the Navier-Stokes system have led to theories of stochastic observables. Importantly, it was found that complete probability density distributions of velocities are governed by the Lundgren-Monin-Novikov (LMN) hierarchy \cite{Lund1967,Monin1967,Novikov1968} that conceptually contains all statistical information.

While it is desirable to have an analogous hierarchy for the scalar model problem, to our best knowledge, such formulation has not been established. Denote by $f^{(N)}= f^{(N)}(t,x_1,v_1,\cdots,x_N,v_N)\ge 0$ and $F^{(N)}=F^{(N)}(t,x_1,v_1,\cdots,x_N,v_N)$ the associated $N$-point probability density function (PDF) and the cumulative density function (CDF) at the points $\{x_k\}_1^N$, respectively, of (\ref{eq:scalar}). That is,
\[
\begin{aligned}
    \text{Prob}(\{u(t,x_k)\in Q^k &\subset \mathbb R | k=1,\cdots, N\}) \\
    &= \int_{Q^N}\cdots \int_{Q^1} f^{(N)}(t,x_1,\tilde v_1,\dots,x_N,\tilde v_N) d\tilde v_1 \cdots d\tilde v_N,
\end{aligned}
\]
and $F^{(N)} = \text{Prob}(\{u(t,x_k)\le v_k | k=1,\dots,N\})$. Of particular interest are the one-point PDF $f=f(t,x,v):=f^{(1)}(t,x_1,v_1)$ and CDF $F=F(t,x,v):=F^{(1)}(t,x_1,v_1)$. Existing approaches to determine $f^{(N)}$ and $F^{(N)}$ include the Monte-Carlo method and the stochastic Galerkin scheme \cite{Xiu2002}. But the evolution of PDFs, which is critical for understanding the decay and asymptotic properties of PDFs, cannot be explicitly characterized by these methods.

Instead, the governing equations for PDFs can be derived by performing ensemble averages to (\ref{eq:scalar}). For scalar problems, master equations for $f$, $F$ and $f^{(2)}$ were reported in previous works \cite{Cho2014,E2000,Pol1995,Tarta2017,Ven2012}. We particularly mention that the analysis of the master equations characterizes the tail probabilities for the velocity gradient and settled a long-standing debate on the stochastic Burgers equation (with $g(u)=u^2/2$) \cite{E2000}.
However, these works either focus on inviscid cases\cite{Tarta2017,Ven2012} or, for the viscous case, result in governing equations with \textit{unknown} forms of conditional averages entering the viscous terms. For instance, the governing equation for $f$ reads as:
\begin{equation} \label{eq:f1}
    f_t + \nabla_x\cdot \left( g'(v) \int_{-\infty}^v f(t,x,\tilde v)d\tilde v \right)_v + \left( \mathbb E_\xi \left[ \epsilon \Delta_x u | u=v \right] f \right)_v = 0,
\end{equation}
which is a linear equation for $f$. For the Burgers equation, this unclosed term $\mathbb E_\xi \left[ \epsilon \Delta_x u | u=v \right]$ can be expressed, in the inviscid limit ($\epsilon\to 0$), in terms of shock-related quantities \cite{E2000}.
Under general cases (with small but finite $\epsilon$), there is no systematic understanding of this term.

The goal of this contribution is to extend (\ref{eq:f1}) to hierarchies of master equations for PDFs and/or CDFs. We develop two approaches. The first (see Section \ref{sec:lmn}) is inspired by the LMN hierarchy that relates multi-point PDFs. The main difference is that here (\ref{eq:scalar}) is a scalar analogue of compressible flows, while the LMN hierarchy was developed based on incompressible Navier-Stokes system. Section \ref{sec:spatial} details a second kind of hierarchy that relates the single-point PDFs of higher-order spatial derivatives. Possible closures based on the hierarchies are discussed in Section \ref{sec:discuss}.

\section{The full hierarchy with $N$-point PDFs} \label{sec:lmn}
In this section, we first write out the master equations for multi-point PDFs $f^{(N)}$ as
\begin{equation} \label{eq:fN}
\begin{split}
    f^{(N)}_t &+ \sum_{k=1}^N \nabla_{x_k}\cdot \left( g'(v_k) \int_{-\infty}^{v_k}f^{(N)}(\cdot,x_k,\tilde v_k) \, d\tilde v_k \right)_{v_k} \\
    &+ \epsilon\sum_{k=1}^N \left( \lim_{x'\to x_k} \Delta_{x'} \int_{\mathbb R} v' f^{(N+1)}(\cdot,x',v') \, dv' \right)_{v_k} = 0.
\end{split}
\end{equation}
It is a linear `kinetic' equation for $f^{(N)}$, with nonlocal effects entering both the advection and viscous terms. Moreover, the viscous term involves an $(N+1)$-point PDF, thus forming an infinite-dimensional hierarchy of equations. 
This new formulation (\ref{eq:fN}) can be viewed as a scalar analogue of the LMN hierarchy for the incompressible Navier-Stokes equations. The main difference is that no divergence-free conditions are generally satisfied in the scalar case. 
For the derivation of (\ref{eq:fN}), the spatial advection term results from the steps in \cite{Pol1995,Ven2012}, whereas the new viscous formula will be derived at the end of this section.

The hierarchy has to be completed by proper initial/boundary data and to obey certain side conditions. For instance, $f^{(N)}$ should be nonnegative, have the normalization and reduction properties
\[
    \int_\mathbb R f(\cdot,v)dv = 1, \quad
    \int_\mathbb R f^{(N)}(\cdot,x_N,v_N)dv_N = f^{(N-1)},
\]
and show `coincidence' behaviors \cite{Lund1967} like
\[
    \lim_{x_2\to x_1} f^{(2)}(t,x_1,v_1,x_2,v_2) = f(t,x_1,v_1)\delta(v_1-v_2),
\]
where $\delta(v_1-v_2)$ is the Dirac Delta-function.

By the definition of CDFs $F^{(N)}$, we integrate (\ref{eq:fN}) in $v_1,\dots,v_N$, interchange the order of integration and the limit, and use the fact $f^{(N)} = F^{(N)}_{v_1\cdots v_N}:= \partial^N F^{(N)}/\partial v_1 \cdots \partial v_N$ to get
\begin{equation} \label{eq:FN}
\begin{split}
    F^{(N)}_t &+ \sum_{k=1}^N g'(v_k) \cdot \nabla_{x_k}F^{(N)} \\
    &+ \epsilon\sum_{k=1}^N \left( \lim_{x'\to x_k} \Delta_{x'} \int_{\mathbb R} v' F^{(N+1)}_{v_k v'}(\cdot,x',v') \, dv' \right) = 0.
\end{split}
\end{equation}
Clearly, this is again a linear hierarchy for $F^{(N)}$. It is remarkable that, unlike (\ref{eq:fN}), the advection term in (\ref{eq:FN}) does not contain nonlocal effects. We mention that for $N=1$, a `method of distribution' has been developed \cite{boso2020,Tarta2017} for the inviscid scalar conservation laws, yielding $F_t+g'(v)\cdot \nabla_x F = \mathcal M$ with an unknown source term $\mathcal M$. By contrast, here we resort to a viscous balance law (\ref{eq:scalar}) so that, with smooth solutions for given $u_0$ and $\epsilon>0$, the master equations (\ref{eq:fN}) and (\ref{eq:FN}) reveal the hierarchical structures.

\begin{remark}
The hierarchical structure is the consequence of the viscous term in (\ref{eq:scalar}). If the scalar balance law is featured with a source term $S=S(u,\nabla_x u)$ that is only dependent on $u$ and $\nabla_x u$: $u_t + g'(u)\cdot \nabla_x u = S$, and assuming smooth solutions $u(t,x)$ exist for $(t,x)\in[0,T]\times D$, then a closed master equation for the joint PDF $q=q(t,x,u,\nabla_x u)$ can be derived. The spatially one-dimensional case is studied in \cite{Ven2012}, and there is no essential difficulty to formulate multidimensional governing equations for $q$. Furthermore, if $S=S(u)$ only relies on $u$, then there exist closed governing equations for $f^{(N)}$ by assuming the existence of smooth solutions of $u$.
\end{remark}

We close this section with a formal derivation of the the hierarchical viscous term in (\ref{eq:fN}), starting from the commonly-considered conditional expectation in (\ref{eq:f1}). We work on the case for $N=1$, but the extension to the general case of $N$ is straightforward. Namely, we need to show
\begin{equation} \label{eq:hier_proof}
    \mathbb E_\xi [\Delta_x u|u=v]f = \lim_{x'\to x}\Delta_{x'} \int_\mathbb R v' f^{(2)}(\cdot,x',v')\, dv'.
\end{equation}
First, it is noticed that the joint PDF of $\Delta_x u$ and $u$, denoted $r^{[\Delta_x u,u]}(a, v)$, can be expressed as (the $t,x$-dependence is omitted)
\[
    r^{[\Delta_x u,u]}(a, v) = \mathbb E_\xi [\delta(\Delta_x u(\xi)-a)\delta(u(\xi)-v)].
\]
Then, by the definition of the conditional expectation, we have
\[
\begin{aligned}
    \mathbb E_\xi [\Delta_x u|u=v]f &= \int ar^{[\Delta_x u,u]}(a, v)\, da
    = \mathbb E_\xi [\Delta_x u(\xi) \delta(u(\xi)-v)] \\
    &= \lim_{x'\to x} \Delta_{x'} \mathbb E_\xi [u(t,x',\xi)\delta(u(t,x,\xi)-v)],
\end{aligned}
\]
where the second line holds by assuming the continuity in $x$ of the conditional expectation. 
Further noting that
\[
\begin{aligned}
    \int_\mathbb R v' &\mathbb E_\xi [\delta(u(t,x,\xi)-v)\delta(u(t,x',\xi)-v')]\, dv' \\
    &= \mathbb E_\xi \left[\delta(u(t,x,\xi)-v) \int_\mathbb R v' \delta(u(t,x',\xi)-v')\, dv' \right] \\
    &= \mathbb E_\xi [u(t,x',\xi)\delta(u(t,x,\xi)-v)]
\end{aligned}
\]
and combining the previous expression, we hence obtain (\ref{eq:hier_proof}) because the two-point PDF $f^{(2)}$ can be expressed as
\[
    f^{(2)}(t,x,v,x',v') = \mathbb E_\xi [\delta(u(t,x,\xi)-v)\delta(u(t,x',\xi)-v')].
\]

\section{New hierarchy with higher-order spatial derivatives} \label{sec:spatial}

There are different approaches to exploit the statistical information of multiple positions. Based on (\ref{eq:scalar}), here we reveal a new set of hierarchical equations for the joint PDFs of higher-order spatial derivatives at a single point. For the sake of simply, assume $d=1$ (the spatially one-dimensional case). The viscous balance law permits smooth solution $u$ and hence the spatial derivatives exist in the classical sense for the sampled initial function $u_0$. Denote $q^{(N)}=q^{(N)}(t,x,a_0,\cdots,a_N)$ as the joint PDFs of $\{\partial_x^k u(t,x)\}_{k=0}^N$, i.e.
\[
\begin{aligned}
    \text{Prob}(\{\partial_x^ku(t,x)\in R^k &\subset \mathbb R | k=0,1,\cdots, N\}) \\
    &= \int_{R^N}\cdots \int_{R^0} q^{(N)}(t,x, \tilde a_0, \dots, \tilde a_N) d\tilde a_0 \cdots d\tilde a_N.
\end{aligned}
\]
The $q^{(N)}$'s are nonnegative and satisfy the normalization and reduction properties like
\[
    \int_\mathbb R q^{(0)}(\cdot,a_0)da_0 = 1, \quad
    \int_\mathbb R q^{(N)}(\cdot,a_N)da_N = q^{(N-1)}.
\]

Then, the master equations for $q^{(N)}$ read as
\begin{equation} \label{eq:qN}
\begin{split}
    q^{(N)}_t =& \sum_{k=0}^N \left( C_k q^{(N)} \right)_{a_k} - \epsilon \sum_{k=0}^{N-2}a_{k+2}q^{(N)}_{a_k} 
    + \left( g'(a_0)A^{(N)} \right)_{a_N} - \epsilon A^{(N)}_{a_{N-1}} - \epsilon B^{(N)}_{a_N},
\end{split}
\end{equation}
where
\[
\begin{aligned}
    A^{(N)} &= A^{(N)}(t,x,a_0,\cdots,a_N) \\
    &= -\int_{-\infty}^{a_N}\left( q^{(N)}_x + \sum_{k=0}^{N-2}a_{k+1}q^{(N)}_{a_k} + \tilde a_N q^{(N)}_{a_{N-1}} \right)(\cdot,\tilde a_N) \, d\tilde a_N,
\end{aligned}
\]
\[
    B^{(N)}=B^{(N)}(t,x,a_0,\cdots,a_N) = \int_\mathbb R A^{(N+1)}(\cdot,a_{N+1})\, da_{N+1},
\]
and the coefficients $C_k$ are defined as
\[
\begin{aligned}
    C_k &= \left( \partial_x^{k+1} g(u) \right)\Big |_{\{\partial_x^i u \gets a_i\}_{i=0}^{k+1}}, \quad k=0,\cdots,N-1, \\
    C_N &= \left( \partial_x^{N+1} g(u) - g'(u)\partial_x^{N+1}u \right)\Big |_{\{\partial_x^i u \gets a_i\}_{i=0}^N}.
\end{aligned}
\]
Here the definition of $C_N$ only contains spatial derivatives of $u$ up to the $N$th order.
Taking $N=2$ as an example, we have $C_0=g'(a_0)a_1$, $C_1=g''(a_0)a_1^2+g'(a_0)a_2$, and $C_2=g'''(a_0)a_1^3+3g''(a_0)a_1a_2$. Clearly, if the flux $g(u)$ is a polynomial, these $C_k$'s can be simplified to some extent.

The governing equations (\ref{eq:qN}) constitute a hierarchy since $B^{(N)}$ contains $A^{(N+1)}$ and hence $q^{(N+1)}$ nonlocally. For clarity, we present below the governing equation for $q^{(0)}$: 
\[
    q_t^{(0)} + \left( g'(a_0) \int_{-\infty}^{a_0} q_x^{(0)}(\cdot,\tilde a_0) d\tilde a_0 \right)_{a_0} = \epsilon \left(\int_{\mathbb{R}} \int_{-\infty}^{a_1} \left( q_x^{(1)} + \tilde{a}_1 q_{a_0}^{(1)} \right)(\cdot,\tilde a_1) d\tilde{a}_1 da_1 \right)_{a_0}.
\]
Since $q^{(0)}(\cdot,a_0)=f(\cdot,a_0)$, comparing the above expression with (\ref{eq:fN}) (with $N=1$), it is evident that the hierarchy of the viscous terms is constructed differently. In (\ref{eq:fN}), there is a limiting process in physical space, whereas the above expression involves a double integral.

\begin{remark}
This method can be extended to higher spatial dimensions ($d>1$) and leads to more complex hierarchies for joint PDFs of the partial derivatives in each spatial direction. Moreover, the idea can also be applied to the Navier-Stokes systems to construct novel hierarchies for fluids.
\end{remark}

To formally derive (\ref{eq:qN}), the main step is to obtain, from the definition
\[
    q^{(N)} = \mathbb E_\xi \left[\delta(u(t,x,\xi)-a_0)\cdots \delta(\partial_x^N u(t,x,\xi) -a_N) \right]
\]
that (the $t,x,\xi$-dependence is omitted)
\begin{equation} \label{eq:chain}
    q^{(N)}_\eta = -\sum_{k=0}^N \partial_{a_k} \mathbb E_\xi \left[ (\partial_\eta \partial_x^k u) \delta(u-a_0)\cdots\delta(\partial_x^N u-a_N) \right], \quad \eta=t,x.
\end{equation}
To prove (\ref{eq:chain}), we take a compactly-supported smooth test function $\varphi=\varphi(a_0,\cdots,a_N)$ and get that  (denoting $da:=da_0\cdots da_N$)
\[
    \int \varphi(a_0,\cdots,a_N)q^{(N)}\, da = \mathbb E_\xi \left[ \varphi(u,\partial_xu,\cdots,\partial_x^Nu) \right],
\]
which implies
\[
\begin{aligned}
   \int \varphi q^{(N)}_\eta da &= \mathbb E_\xi \left[ \sum_{k=0}^N \varphi_{\partial_x^k u} \partial_\eta \partial_x^ku \right]
   = \mathbb E_\xi \left[ \sum_{k=0}^N \left( \int \varphi_{a_k} \prod_{j=0}^N\delta(\partial_x^j u-a_j)da \right) \partial_\eta \partial_x^ku \right] \\
   &= \int \sum_{k=0}^N \varphi_{a_k} \mathbb E_\xi \left[ (\partial_\eta \partial_x^ku) \prod_{j=0}^N\delta(\partial_x^j u-a_j) \right]da
   = \int \varphi R da.
\end{aligned}
\]
Here $R=R(t,x,a_0,\cdots,a_N)$ is the right-hand side of (\ref{eq:chain}). Since $\varphi$ is arbitrary, (\ref{eq:chain}) follows immediately. With similar arguments it can be shown that
\begin{equation} \label{eq:replace}
    \mathbb E_\xi [h(u,\dots,\partial_x^N u) \delta(u-a_0)\cdots \delta(\partial_x^N u-a_N)] = h(a_0,\dots,a_N) q^{(N)}
\end{equation}
holds for any function $h$ that only relies on $u,\dots,\partial_x^N u$.

Then one derives (\ref{eq:qN}) by taking $\eta:=t$ in (\ref{eq:chain}), substituting (\ref{eq:scalar}) for $u_t$, and then using (\ref{eq:replace}). After these steps, the terms that still require handling include forms $\mathbb E_\xi[(\partial_x^k u) \delta(u-a_0)\cdots\delta(\partial_x^Nu-a_N)]$ with $k=N+1,\ N+2$. Indeed, the two forms are denoted as $A^{(N)}$ and $B^{(N)}$ in (\ref{eq:qN}), whose expressions therein are again the consequence of (\ref{eq:chain}) (with $\eta:=x$) and (\ref{eq:replace}).

\section{Discussions} \label{sec:discuss}

The spatial dependency rooted in the viscous term in (\ref{eq:scalar}) prevent the derivation of a closed governing equation for the PDF. In Sections \ref{sec:lmn} and \ref{sec:spatial}, two hierarchical construction approaches are developed, involving multi-point PDFs/CDFs (\ref{eq:fN})(\ref{eq:FN}) and single-point spatial-derivative PDFs (\ref{eq:qN}). These approaches suggest a flexible framework for capturing the statistical information of correlations across different spatial locations, potentially offering better adaptability to measurement methods in real-world scenarios.

The practical implication of these hierarchies also lies in providing a pathway for the development of closures, which requires extra assumptions on higher-order PDFs. This resembles in spirit the celebrated BBGKY hierarchy of kinetic gas theory, as noted in \cite{Lund1967}. For the one-point PDF $f$ and CDF $F$, if we assume spatial independency (denoting $w_i=(x_i,v_i)$)
\[
    f^{(2)}(t,w_1,w_2) = f(t,w_1)f(t,w_2) \iff F^{(2)}(t,w_1,w_2) = F(t,w_1)F(t,w_2),
\]
then a closed equation for $F$ is derived based on (\ref{eq:FN}) as
\begin{equation} \label{eq:F_ind}
    F_t+g'(v)\cdot \nabla_xF + \epsilon \left( \Delta_x \int_\mathbb R \tilde vF_v(t,x,\tilde v)\, d\tilde v \right) F_v = 0,
\end{equation}
where the continuity of $\Delta_x \int_\mathbb R \tilde vF_v(t,x,\tilde v)d\tilde v$ in $x$ has been assumed.
However, this Ansatz may be oversimplified in most cases. Another approach originating from physical investigation takes the approximation \cite{Lund1967}
\[
    f^{(3)}(t,w_1,w_2,w_3) = \prod_{k=1,2,3} f(t,w_k) + \sum_{i,j,k} f(t,w_i)\tilde f^{(2)}(t,w_j,w_k),
\]
where the indices $1\le i,j,k\le 3$ are mutually distinct and
\[
    \tilde f^{(2)}(t,w_1,w_2) = f^{(2)}(t,w_1,w_2) - f(t,w_1)f(t,w_2).
\]
Obviously, this Ansatz results in closed governing equations for $f$ and $f^{(2)}$ out of (\ref{eq:fN}), or equivalently, for $F$ and $F^{(2)}$ out of (\ref{eq:FN}).

Future work can be directed to more thoroughly explore closure strategies based on the hierarchies developed here (Eqs. (\ref{eq:fN}), (\ref{eq:FN}), and (\ref{eq:qN})). Theoretical analyses of the resulting closed systems (e.g., (\ref{eq:F_ind})), as well as numerical treatments, are also promising directions for future research.

\section*{Acknowledgements}
Funding by the Deutsche Forschungsgemeinschaft (DFG, German Research Foundation) - SPP 2410 \textit{Hyperbolic Balance Laws in Fluid Mechanics: Complexity, Scales, Randomness (CoScaRa)} is gratefully acknowledged.

\bibliographystyle{amsplain}
\bibliography{references}

\end{document}